\documentclass{amsart}
\usepackage{amsfonts,amscd}

\newtheorem{theorem}{Theorem}[section]
\newtheorem{lemma}[theorem]{Lemma}
\newtheorem{corollary}[theorem]{Corollary}
\newtheorem{proposition}[theorem]{Proposition}
\theoremstyle{remark}

\theoremstyle{definition}

\numberwithin{equation}{section}
\makeatother

\begin{document}

\title[Szeg\"{o}'s theorem and outers for
noncommutative
$H^p$]{Applications
of the Fuglede-Kadison determinant: \\
 Szeg\"{o}'s theorem and outers\\ for noncommutative $H^p$}

\date{\today}

\author{David P. Blecher}
\address{Department of Mathematics, University of Houston, Houston, TX
77204-3008}
\email[David P. Blecher]{dblecher@math.uh.edu}
 \author{Louis E. Labuschagne}
\address{Department of Mathematical Sciences, P.O. Box
 392, 0003 UNISA, South Africa}
\email{labusle@unisa.ac.za}
\thanks{*Blecher was partially supported by grant DMS 0400731
from the National Science Foundation, and Labuschagne
partially by a National Research Foundation Focus Area Grant}

\begin{abstract}
 We first use properties of the Fuglede-Kadison
 determinant on $L^p(M)$, for a finite von Neumann algebra $M$, to
  give several useful variants of
the noncommutative Szeg\"{o} theorem for $L^p(M)$, including the one
usually attributed to Kolmogorov and Krein. As an application, we
solve the longstanding open problem concerning the noncommutative
generalization, to Arveson's noncommutative $H^p$ spaces, of the  
famous `outer factorization' of
functions $f$ with $\log |f|$ integrable.   Using the
Fuglede-Kadison
 determinant, we also generalize
many other classical results concerning outer functions.
\end{abstract}

\maketitle

\section{Introduction}

It has long been of great importance to operator theorists and
operator
algebraists to find noncommutative analogues of the classical
`inner-outer factorization' of analytic functions.
 We recall some classical results:
If $f \in L^1$
with $f \geq 0$, then $\int  \log f \,
> - \infty$ iff $f = |h|$ for an outer $h \in H^1$ (iff $f = |h|^p$
for an outer $h \in H^p$).   We will call this
the {\em Riesz-Szeg\"{o} theorem} since it is often
attributed to one or other of these authors.   If $f \in L^1$ with
$\int  \log |f| \,
> - \infty$, then $f = uh$, where $u$ is  unimodular
 and $h$ is outer.  Moreover, outer functions may be {\em defined} in
terms of a simple equation involving $\int  \log |f|$.   Such
results are usually treated as consequences of the classical {\em
Szeg\"o theorem}, which is really a distance formula in terms of the
entropy $\exp(\int  \log |f|)$, and which in turn is intimately
related to the {\em Jensen inequality} (see e.g.\ \cite{Hobk}). In
the noncommutative situation one wishes, for example,
 to find conditions on a positive operator $T$
which imply that $T = |S|$ for an operator $S$ which is in a
`noncommutative Hardy class',
or, even better, which is `outer' in some sense.
There are too many
such results in the literature to attempt a listing of them
(see e.g.\ \cite[p.\ 1495]{PX}).
 Indeed this is an active and
important research area which has links to many exciting parts of
mathematics. Interestingly, central parts of this topic still seem
to be poorly understood. As an example of this, we cite the main and
now classical result of \cite{Dev},  concerning a Riesz-Szeg\"o like
factorization of a class of $B(H)$-valued functions on the unit interval, which
has resisted generalization in some important directions. In our
paper we generalize the classical results above to the noncommutative
$H^p$ spaces associated with Arveson's remarkable {\em subdiagonal
algebras} \cite{AIOA}.  Our generalization solves an old open
problem (see the discussion in \cite[lines 8-12, p.\ 1497]{PX}, and
\cite[p.\ 386]{MMS}).  The approach which we take has been
unavailable until now (since it relies ultimately on the recent
solution in \cite{LL3} of a 40 year old open problem from
\cite{AIOA}).  Moreover, the approach is very faithful to the
original classical function theoretic route (see e.g.\ \cite{Hobk}),
proceeding via noncommutative Szeg\"o theorems.

 In the last several years we have attempted to demonstrate
that all the results in (a particular survey \cite{SW} of) the
`generalized $H^p$-theory' for abstract function algebras from the
1960s, extend in  an extremely complete and literal fashion, to the
noncommutative setting of Arveson's subdiagonal subalgebras of von
Neumann algebras \cite{AIOA}.  This may be viewed as a very natural
merging of generalized Hardy space, von Neumann algebra, and
noncommutative
$L^p$ space, techniques.   See our recent survey \cite{BL5} for an
overview of this work. The present paper completes the
noncommutative extension of the basic Hardy space theory. As posited
by Arveson, one should use the Fuglede-Kadison determinant
$\Delta(a) = \exp(\tau(\log |a|))$ where $\tau$ is a trace, as a
natural replacement in the noncommutative case for the quantity
$\int \log f$ above.  We use properties of the Fuglede-Kadison
determinant to give several useful variants of the noncommutative
Szeg\"{o} theorem for $L^p(M)$, including the one usually attributed
to Kolmogorov and Krein. As applications, we generalize the
noncommutative Jensen inequality, and generalize many of the
classical results concerning outer functions, to the noncommutative
$H^p$ context.  Some of these will be described in more detail later
in this introduction.

We now review some of the definitions and notation, although we
strongly suggest that the reader glance though our survey article
\cite{BL5} first for background, motivation, history, etc.
 For a set ${\mathcal S}$, we write ${\mathcal
S}_+$ for the set $\{ x \in {\mathcal S} : x \geq 0 \}$.  We assume
throughout that $M$ is a von Neumann algebra possessing a faithful
normal tracial state $\tau$. The existence of such $\tau$ implies
that $M$ is a so-called {\em finite} von Neumann algebra, and that
if $x^* x = 1$ in $M$, then $x x^* = 1$ too. Indeed, for any $a, b
\in M$, $ab$ will be invertible precisely when $a$ and $b$ are
separately invertible. We will also need to use a well known fact
about inverses of an unbounded
operator $T$, and in our case $T$ will be positive, selfadjoint,
closed, and densely defined.
We recall that $T$ is {\em
bounded below} if for some $\lambda > 0$ one has $\|T(\eta)\| \geq
\lambda\|\eta\|$ for all $\eta \in \mathrm{dom}(T)$.
This is equivalent to demanding that
$|T| \geq \epsilon 1$ for some $\epsilon > 0$, and of course in this
case, $|T|$ has a bounded positive inverse.

 A (finite maximal) {\em subdiagonal subalgebra} of $M$ is a weak*
 closed unital subalgebra $A$ of
$M$ such that if $\Phi$ is the unique conditional expectation
guaranteed by \cite[p.\ 332]{Tak} from $M$ onto $A \cap A^*
\overset{def}{=} {\mathcal D}$ which is trace preserving (that is,
$\tau \circ \Phi = \tau$), then:
\begin{equation} \label{Eq2}
 \Phi(a_1 a_2) = \Phi(a_1) \,  \Phi(a_2) , \; \; \; a_1, a_2 \in A .
\end{equation}
One also must impose one further condition on $A$.  There is a
choice of at least eight equivalent, but quite different looking,
conditions \cite{BL5}; Arveson's original one (see also \cite{E}) is
that $A + A^*$ is weak* dense in $M$.  In the classical function
algebra setting \cite{SW}, one assumes that ${\mathcal D} = A \cap
A^*$ is one dimensional, which forces $\Phi = \tau(\cdot) 1$. If in
our setting this is the case, then we say that $A$ is {\em
antisymmetric}.

The simplest example of a maximal subdiagonal algebra is the upper
triangular matrices $A$ in $M_n$. Here  $\Phi$ is the expectation
onto the main diagonal. There are much more interesting examples
from free group von Neumann algebras, Tomita-Takesaki theory, etc
(see e.g.\ \cite{AIOA,Zs,MMS}).  On the other end of the spectrum,
$M$ itself is a maximal subdiagonal algebra (take $\Phi = Id$).  It
is therefore remarkable that so much of the classical $H^p$ theory
does extend to {\em all} maximal subdiagonal algebras. However the
reader should not be surprised to find
one or two results below which do impose
restrictions on the size of ${\mathcal D}$.

 By analogy with the classical case, we set $A_0 = A \cap
{\rm Ker}(\Phi)$ and  set $H^p$ or $H^p(A)$ to be $[A]_p$, the
closure of $A$ in the noncommutative $L^p$ space $L^p(M)$, for $p
\geq 1$.   More generally we write $[{\mathcal S}]_p$ for this
closure of any subset ${\mathcal S}$.  We will often view $L^p(M)$
inside $\widetilde{M}$, the set of unbounded, but closed and densely
defined, operators on $H$ which are affiliated to $M$.
In the present context this is a
$*$-algebra with respect to the `strong' sum and product (see
Theorem I.28 and the example following it in \cite{Terp}).
We order $\widetilde{M}$ by its cone of positive (selfadjoint)
 elements.
The trace
$\tau$ extends naturally to the
positive operators in
$\widetilde{M}$.  If $1 \leq p < \infty$, then $L^p(M,\tau) =
 \{a \in \widetilde{M} :  \tau(|a|^p) < \infty\}$,
equipped with the norm $\|\cdot\|_p = \tau(|\cdot|^p)^{1/p}$ (see
e.g.\ \cite{Nel,FK,Terp,PX}).  We abbreviate $L^p(M,\tau)$ to
$L^p(M)$.

Arveson's Szeg\"o formula is:
$$\Delta(h) = \inf \{ \tau(h |a+d|^2)
: a \in A_0 , d \in {\mathcal D} , \Delta(d) \geq 1 \}$$ for all $h
\in L^1(M)_+$.  Here $\Delta$ is the {\em Fuglede-Kadison
determinant}, originally defined on $M$ by  $\Delta(a) = \exp
\tau(\log |a|)$ if $|a| > 0$, and otherwise, $\Delta(a) = \inf \,
\Delta(|a| + \epsilon 1)$, the infimum taken over all scalars
$\epsilon > 0$ (see \cite{FKa,AIOA}).  We will
discuss this determinant in more detail in Section 2.   Unfortunately,
the just-stated noncommutative
Szeg\"o formula, and the (no doubt  more important)
associated {\em Jensen's inequality}
$$\Delta(\Phi(a)) \leq \Delta(a) , \qquad a \in A ,$$
resisted proof for nearly 40 years, although Arveson did prove them
in his extraordinary original paper \cite{AIOA} for the examples
that he was most interested in. In 2004, the second author proved in
\cite{LL3} that all maximal subdiagonal algebras satisfy these
formulae.  Settling this old open problem opened up the theory, and
in particular enabled the present work.

An element $h \in H^p$ is said to be {\em outer} if $1 \in [h A]_p$.
This definition is in line with e.g.\ Helson's definition
of outers in the matrix valued case he considers in
\cite{Hel}.   We now state a sample of our results about outers.
 For example, we are able to
improve on the factorization
theorems from e.g.\ \cite[Section 3]{BL3} in several ways:
namely we show that if $f \in L^p(M)$ with $\Delta(f) > 0$ then
$f$ may be essentially uniquely factored $f = uh$ with
$u$ unitary and $h$ outer.  There is a much more obvious converse
to this, too.
Moreover we now have an explicit formula
for the $u$ and $h$.    We refer to a factorization
$f = uh$  of this form as a {\em Beurling-Nevanlinna factorization}.
It follows that in this case if $f \geq 0$ then
$f = |h|$ with $h$ outer.
  This gives a
solution to the problem posed in \cite[Remark after Theorem
8.1]{PX}, and in \cite[p.\ 386]{MMS}.  
If $h \in H^p$,
and $h$ is outer then $\Delta(h) = \Delta(\Phi(h))$.
Moreover, a converse is true: if  
 $\Delta(h) = \Delta(\Phi(h)) > 0$ then $h$ is outer.
It follows that under some restrictions on ${\mathcal D} =
A \cap A^*$, $h$ is outer iff $\Delta(h) = \Delta(\Phi(h)) > 0$.

Some historical remarks: there are many factorization theorems for
subdiagonal algebras in the literature (see e.g.\
\cite{AIOA,MMS,Sai,MW,PX}), but as far as we know there are no
noncommutative factorization results\footnote{With the exception of
the factorization results in our previous paper \cite{BL3}, which
are in many senses improved upon here.} involving outers or the
Fuglede-Kadison determinant.  We mention for example Arveson's
original factorization result from \cite{AIOA}, or Marsalli and
West's Riesz factorization of any $f \in H^1$ as a product $f = g h$
with $g \in H^p, h \in H^q$, $\frac{1}{p} + \frac{1}{q} = 1$. Some
also require rather stronger hypotheses, such as $f^{-1} \in L^2(M)$
(see e.g.\ \cite{MMS}).

Another important historical remark,
is that the commutative case of most of the topics in
our paper was settled in \cite{Nak}.  While this paper certainly gave
us motivation to persevere in our endeavor, we follow completely different
lines, and indeed the results work out
rather differently too.  In particular, the quantity $\tau(\exp(\Phi(\log|f|)))$,
which plays a central role in most of the results in \cite{Nak},
seems to us to be unrelated to outers or factorization in the noncommutative
setting.  In passing, we remark that numerical experiments
do seem to confirm the existence of a Jensen inequality
 $\tau(\exp(\Phi(\log|a|))) \geq \tau(|\Phi(a)|)$ for subdiagonal
algebras.  However, even if correct, it is not clear
how this might be useful to the present work.   We therefore
defer such considerations to a future investigation.

Finally, we remark that there are many other, more recent,
generalizations of $H^\infty$, based around multivariable analogues
of the Sz-Nagy-Foia\c s model theory for contractions. In essence,
the unilateral shift is replaced by left creation operators on some
variant of Fock space.  Many prominent researchers are currently
intensively pursuing these topics, they  are very important and are
evolving in many directions. Although these theories also contain
variants of  Hardy space theory, they are quite far removed from
subdiagonal algebras. For example, if one compares Popescu's theorem
of Szeg\"o type from \cite[Theorem 1.3]{Pop} with the Szeg\"o
theorem for subdiagonal algebras discussed here, one sees that they
are only related in a very formal sense.

\section{Properties of the Fuglede-Kadison determinant}

The Fuglede-Kadison determinant $\Delta$, and its amazing properties, is
perhaps the main tool in the noncommutative $H^p$ theory.
In \cite{FK}, Fuglede
and Kadison study the determinant as a function on $M$.  In
the next paragraph we will define it for elements of $L^q(M)$ for any $q > 0$.
In fact, as was pointed out to us by Quanhua Xu, L. G. Brown
investigated the determinant and its properties in the early 1980s,
on a much larger class than $L^q(M)$ (see
\cite{Brme}); indeed recently Haagerup and Schultz
have thoroughly explicated the basic theory of this determinant
for a very general class of $\tau$-measurable operators (see \cite{HS})
as part of Haagerup's amazing attack on the invariant subspace problem relative
to a finite von Neumann algebra.

For our purposes, we will define the Fuglede-Kadison determinant for an element
$h \in L^q(M)$,
for any $q > 0$, as follows.
We set $\Delta(h) = \exp\tau(\log |h|)$
if $|h| > \epsilon 1$ for some $\epsilon > 0$, and otherwise, $\Delta(h) =
\inf \, \Delta(|h| + \epsilon 1)$, the infimum taken over all scalars
$\epsilon > 0$.  To see that this is well-defined,
 we adapt the argument in the
third paragraph of \cite[Section 2]{BL3},
making use of the Borel
functional calculus for unbounded operators applied to the inequality
$$0 \leq \log t \leq \frac{1}{q}t^q , \qquad t \in [1,\infty).$$
Notice that for any $0 < \epsilon < 1$, the function $\log t$ is
bounded on $[\epsilon,1]$. So given $h \in L^1(M)_+$ with $h \geq
\epsilon 1$, it follows that $(\log h) e_{[0,1]}$ is similarly
bounded. Moreover the previous centered equation ensures that $0
\leq (\log h) e_{[1,\infty)} \leq \frac{1}{q}h^q e_{[1,\infty)} \leq
\frac{1}{q}h^q.$ Here $e_{[0,\lambda]}$ denotes the spectral
resolution of $h$. Thus if $h \in L^q(M)$ and $h \geq \epsilon$ then
$\log h \in L^1(M)$.

The following are the basic properties
of this extended determinant which we shall need.  Full
proofs may be found in \cite{HS}, which are valid for a very
general class of unbounded operators
(see also \cite{BL5} for another (later) proof
for the $L^p(M)$ class).
We will often use these results silently in the next few sections.

\begin{theorem} \label{Fkd}  If $p > 0$ and $h \in L^p(M)$  then
\begin{itemize}
\item [(1)]  $\Delta(h) = \Delta(h^*) = \Delta(|h|)$.
\item [(2)]   If $h \geq g$ in $L^p(M)_+$ then $\Delta(h) \geq \Delta(g)$.
\item [(3)]   If $h \geq 0$ then $\Delta(h^q) = \Delta(h)^q$ for any $q > 0$.
\item [(4)]  $\Delta(h b) = \Delta(h)\Delta(b) = \Delta(b h)$ for any
$b \in L^q(M)$ and any $q > 0$.
\end{itemize}
\end{theorem}

\section{Szeg\"o's formula revisited}

Throughout this section, $A$ is a maximal subdiagonal algebra in
$M$. We consider versions of  Szeg\"o's formula valid in $L^p(M)$
rather than $L^2(M)$.   We will also prove a generalized Jensen
inequality, and show that the classical Verblunsky/Kolmogorov-Krein
strengthening of Szeg\"o's formula extends even to the
noncommutative context.

It is proved in \cite{BL4} that for
 $h \in L^1(M)_+$ and $1 \leq p < \infty$, we have
$$\Delta(h) = \inf \{ \tau(h |a+d|^p) : a \in A_0 , d \in {\mathcal
D} , \Delta(d) \geq 1 \} .$$
We now prove some perhaps more useful variants of this formula.

\begin{lemma} \label{vnc}  If $h \in L^q(M)_+$ and $0 < p, q < \infty$,
we have
$$\Delta(h) = \inf \{ \tau(|h^{\frac{q}{p}} b|^p)^{\frac{1}{q}} : b \in
M_+, \Delta(b) \geq 1 \} = \inf \{ \tau( |b
h^{\frac{q}{p}}|^p)^{\frac{1}{q}} : b \in M_+, \Delta(b) \geq 1 \}
.$$ The infimums are realized on the commutative von Neumann
subalgebra $M_0$ generated by $h$, and are unchanged if in addition
we also require $b$ to be invertible in $B$. \end{lemma}

\begin{proof}  That the two infimums in the displayed equation are
equal follows from the fact that $\Vert x \Vert_p =
\Vert x^* \Vert_p$ for $x \in L^p(M)$ (see \cite{FK}).
 Thus we just prove the first equality in that line.

For
 $b \in M_+, \Delta(b) \geq 1$, we have by Theorem \ref{Fkd} (3) that
$$\Delta(|h^{q/p}b|^p) =\Delta(|h^{q/p}b|)^p = \Delta(h^{q/p}b)^p .$$
Consequently, using facts from Theorem \ref{Fkd} again, we have
$$\tau(|h^{q/p}b|^p)\geq \Delta(|h^{q/p}b|^p)
= [\Delta(h^{q/p})\Delta(b)]^p \geq \Delta(h^{q/p})^p = \Delta(h)^q.$$
To complete the proof, it suffices to find, given $\epsilon > 0$, an
invertible $b$ in $(M_0)_+$, the von Neumann algebra generated by $h$
(see the first paragraph of Section 2), with $\Delta(b) \geq 1$ and
$\tau(|h^{\frac{q}{p}} b|^p)^{\frac{1}{q}} < \Delta(h) + \epsilon$.
But for any $b \in (M_0)_+$ we have $|h^{q/p}b|^p = h^qb^p$ by
commutativity, and then the result follows from an analysis of
Arveson's original definition of $\Delta(h)$ (see (2.1) in
\cite{BL2}). In particular since $\Delta(h^q) = \inf \{ \tau(h^q
b^p) : b \in (M_0)_+, \Delta(b) \geq 1 \}$ by \cite[Theorem
2.1]{BL2}, an application of Theorem \ref{Fkd} (3) ensures that
$\Delta(h) = \inf \{ \tau(h^q b^p)^{\frac{1}{q}} : b \in (M_0)_+,
\Delta(b) \geq 1 \}$.
\end{proof}

\begin{corollary} \label{vnc2}  If $h \in L^q(M)_+$ and $0 < p, q < \infty$,
we have
\begin{eqnarray*}
\Delta(h) &=& \inf \{\tau(|h^{\frac{q}{p}} a|^p)^{\frac{1}{q}}
: a \in A, \Delta(\Phi(a)) \geq 1 \}\\
&=& \inf \{\tau(|a h^{\frac{q}{p}}|^p)^{\frac{1}{q}} : a \in A,
\Delta(\Phi(a)) \geq 1 \}.
\end{eqnarray*}
The infimums are unchanged if we also require $a$ to be invertible
in $A$, or if we require $\Phi(a)$ to be invertible in ${\mathcal
D}$.
\end{corollary}

\begin{proof}   That the two infimums in the displayed equation are
equal follows from the fact that $\Vert x \Vert_p =
 \Vert x^* \Vert_p$ for $x \in L^p(M)$ (see \cite{FK}),
 and by replacing $A$ with $A^*$, which is also subdiagonal.
 Thus we just prove the first equality in that line.

For $a \in A, \Delta(\Phi(a)) \geq
 1$ we have
$$\tau(|h^{\frac{q}{p}} a|^p)^{\frac{1}{q}} =
\tau(|a^* h^{\frac{q}{p}}|^p)^{\frac{1}{q}} =
\tau(||a^*| h^{\frac{q}{p}}|^p)^{\frac{1}{q}} \geq \Delta(h) ,$$ by Lemma
\ref{vnc}, since $\Delta(|a^*|) = \Delta(a^*) = \Delta(a) \geq
\Delta(\Phi(a)) \geq 1$ (using Jensen's inequality). Thus
$\Delta(h)$ is dominated by the first infimum. On the other hand, by
the previous result there is an invertible $b \in M_+$ with
$\Delta(b) \geq 1$ and $\tau(|b h^{\frac{q}{p}}|^p)^{\frac{1}{q}} <
\Delta(h) + \epsilon$. By factorization, we can write $b = |a^*|$
for an invertible $a$ in $A$, and by Jensen's formula
\cite{AIOA,LL3} we have
$\Delta(\Phi(a)) = \Delta(a) = \Delta(a^*) = \Delta(b) \geq 1$. Then
$$\tau(|h^{\frac{q}{p}} a|^p)^{\frac{1}{q}} =
\tau(|a^* h^{\frac{q}{p}}|^p)^{\frac{1}{q}} =
\tau(|b h^{\frac{q}{p}}|^p)^{\frac{1}{q}} < \Delta(h) + \epsilon .$$
We leave the remaining details to the reader.
\end{proof}

\begin{corollary} \label{vncc2}  {\rm (Generalized Jensen inequality)} \
 Let $A$ be a maximal subdiagonal algebra. For any $h \in H^1$
 we have $\Delta(h) \geq \Delta(\Phi(h))$.
  \end{corollary}

\begin{proof}  Using the $L^1$-contractivity of $\Phi$
we get $$\tau(||h|a|) = \tau(|ha|)
\geq \tau(|\Phi(ha)|) = \tau(||\Phi(h)|\Phi(a)|) , \qquad  a \in A .$$
Taking the infimum over such $a$  with $\Delta(\Phi(a)) \geq 1$,
we obtain from Corollary \ref{vnc2}, and Theorem \ref{vnc}
applied to ${\mathcal D}$, that $\Delta(h)
= \Delta(|h|) \geq \Delta(|\Phi(h)|) = \Delta(\Phi(h))$. \end{proof}

 We recall that although $L^p(M)$ is not a normed space
if $1 > p > 0$, it is a so-called linear metric space with metric given by
$\|x - y\|_p^p$ for any $x, y \in L^p$ (see \cite[4.9]{FK}). Thus although
the unit ball may not be convex, continuity still respects all elementary
linear operations.

\begin{corollary} \label{vnc3}  Let $h \in L^q(M)_+$ and
$0 < p, q < \infty$.
If $h^{\frac{q}{p}} \in [h^{\frac{q}{p}} A_0]_p$, then
$\Delta(h) = 0$.  Conversely, if $A$ is antisymmetric
and $\Delta(h) = 0$, then $h^{\frac{q}{p}} \in [h^{\frac{q}{p}} A_0]_p$.
Indeed if $A$ is antisymmetric, then
$$\Delta(h) = \inf \{\tau(|h^{\frac{q}{p}} (1 - a_0)|^p)^{\frac{1}{q}} :
a_0 \in A_0 \} .$$
\end{corollary}

\begin{proof}  The first assertion follows by
taking $a$ in the infimum in Corollary \ref{vnc2} to be
of the form $1 - a_0$ for $a_0 \in A_0$.

If $A$ is antisymmetric, and if $t \geq 1$ with
$\tau(|h^{\frac{q}{p}} (t1 + a_0)|^p)^{\frac{1}{q}} < \Delta(h) +
\epsilon$, then $\tau(|h^{\frac{q}{p}} (1 + a_0/t)|^p)^{\frac{1}{q}}
< \Delta(h) + \epsilon$. From this the last assertion follows that
the infimum's in Corollary \ref{vnc2} can be taken over terms of the
form $1 + a_0$ where $a_0 \in A_0$. If this infimum was $0$ we could
then find a sequence $a_n \in A_0$ with $h^{\frac{q}{p}} (1 + a_n)
\to 0$ with respect to $\Vert \cdot \Vert_p$. Thus $h^{\frac{q}{p}}
\in [h^{\frac{q}{p}}A_0]_p$.
\end{proof}

{\bf Remark.}  The converse in the last result is false for general
maximal subdiagonal algebras (e.g.\ consider $A = M = M_n$).

For a general maximal subdiagonal algebra $A$, and $h \in L^1(M)_+$
we define $\delta(h) = \inf \{\tau(|h^{\frac{1}{2}} (1 - a_0)|^2) :
a_0 \in A_0 \} .$   Note that $\delta(h) = \tau(h)$ in the example
$M = A$.

\medskip

We close this section with the following version of the Szeg\"o
formula valid for general positive linear functionals.  Although the
classical version of this theorem is usually attributed to
Kolmogorov and Krein, we have been informed by Barry Simon that
Verblunsky proved it first, in the mid 1930's (see e.g.\ \cite{Ver}):

\begin{theorem} \label{SKK}  {\rm (Noncommutative
Szeg\"o-Verblunsky-Kolmogorov-Krein theorem)}  \ Let $\omega$ be a positive linear
functional on $M$, and let $\omega_n$ and
$\omega_s$ be its normal and singular parts respectively, with
$\omega_n = \tau(h \, \cdot \, )$ for $h \in
L^1(M)_+$. If $\mathrm{dim}(\mathcal{D}) < \infty$, then
$$\Delta(h) = \inf \{\omega(|a|^2) : a \in A, \Delta(\Phi(a)) \geq 1 \}.$$
The infimum remains unchanged if we also require
$\Phi(a)$ to be invertible in ${\mathcal
D}$.
\end{theorem}

\begin{proof}
Suppose that $\mathrm{dim}(\mathcal{D}) < \infty$. All terminology
and notation will be as in Lemma 3.2 of \cite{BL4}, the preamble to
the proof of the noncommutative F. \& M. Riesz theorem \cite{BL4}.
For the sake of clarity we pause to highlight the most important of these.
If $(\pi_\omega, H_\omega, \Omega_\omega)$ is the GNS representation
engendered by $\omega$, there exists a central
projection $p_0$ in $\pi_{\omega}(M)''$ such that for any $\xi, \psi
\in H_{\omega}$ the functionals $a \mapsto
\langle\pi_{\omega}(a)p_0\xi, \psi\rangle$ and $a \mapsto
\langle\pi_{\omega}(a)(1-p_0)\xi, \psi\rangle$ on $M$ are
respectively the normal and singular parts of the functional $a
\mapsto \langle\pi_{\omega}(a)\xi, \psi\rangle$ \cite[III.2.14]{Tak}. In this
representation $\Omega_0$ will denote the orthogonal projection of
$\Omega_\omega$ onto the closed subspace
$\overline{\pi_\omega(A_0)\Omega_\omega}$.

Let $d \in \mathcal{D}$ be  given. We may of course select a sequence
$(f_n) \subset A_0$ so that $\lim_{n \to \infty}\pi_\omega(f_n)\Omega_\omega
=\Omega_0$. By the ideal property of $A_0$ and continuity, it then follows
that $$\pi_\omega(d)\Omega_0 = \lim_{n \to \infty}\pi_\omega(df_n)\Omega_\omega
\in \overline{\pi_\omega(A_0)\Omega_\omega}.$$

Once again using the ideal property of $A_0$, the fact that
$\Omega_{\omega}-\Omega_0 \perp
\overline{\pi_\omega(A_0)\Omega_\omega}$ now forces
$\langle\pi_{\omega}(d)(\Omega_{\omega}-\Omega_0),
\pi_\omega(a)\Omega_{\omega}\rangle = \langle
\Omega_{\omega}-\Omega_0 , \pi_\omega(d^*a)\Omega_{\omega}\rangle =
0$ for every $a \in A_0$. Therefore
$$\pi_{\omega}(d)(\Omega_{\omega}-\Omega_0) \perp
\overline{\pi_\omega(A_0)\Omega_\omega}.$$  From the facts in the
previous two centered equations, it now follows that
$\pi_{\omega}(d)\Omega_0$ is the orthogonal projection of
$\pi_{\omega}(d)\Omega_{\omega}$ onto
$\overline{\pi_\omega(A_0)\Omega_\omega}$. Using this fact we may
now repeat the arguments of \cite[Lemma 3.2(a) \& (b)(i)]{BL4} for
the functional
$$\omega_d(\cdot) =
\langle\pi_\omega(\cdot)\pi_{\omega}(d)(\Omega_{\omega}-\Omega_0),
\pi_{\omega}(d)(\Omega_{\omega}-\Omega_0)\rangle$$ to conclude that
$\omega_d$ is normal, with $p_0(\pi_{\omega}(d)(\Omega_{\omega}-
\Omega_0)) = \pi_{\omega}(d)(\Omega_{\omega}-\Omega_0)$, where $p_0$
is the central projection in $\pi_\omega(M)''$ mentioned above, and
also: $p_0 \pi_\omega(d) (\Omega_{\omega}-\Omega_0) \perp
p_0 \pi_\omega(A_0)\Omega_\omega$ and $p_0(\pi_{\omega}(d)\Omega_0)$
is the orthogonal projection of $p_0(\pi_{\omega}(d)\Omega_{\omega})$
onto $p_0(\overline{\pi_\omega(A_0)\Omega_\omega})$. Thus we arrive
at the fact that

\begin{eqnarray*}
\inf_{a \in A_0} \omega(|d + a|^2) &=& \inf_{a \in A_0}
\langle\pi_{\omega}(d)\Omega_{\omega}+\pi_\omega(a)\Omega_\omega,
\pi_{\omega}(d)\Omega_{\omega}+\pi_\omega(a)\Omega_\omega\rangle \\
& = & \inf_{a \in A_0} \, \Vert \pi_{\omega}(d)\Omega_{\omega}-
\pi_\omega(a)\Omega_\omega \Vert^2  \\ &=&
\langle\pi_{\omega}(d)(\Omega_{\omega}-\Omega_0),
\pi_{\omega}(d)(\Omega_{\omega}-\Omega_0)\rangle \\
&=& \langle p_0\pi_{\omega}(d)(\Omega_{\omega}-\Omega_0),
p_0\pi_{\omega}(d)(\Omega_{\omega}-\Omega_0)\rangle \\
&=& \inf_{a \in A_0}
\langle p_0 \pi_{\omega}(d)\Omega_{\omega}+p_0\pi_\omega(a)\Omega_\omega,
p_0\pi_{\omega}(d)\Omega_{\omega}+p_0\pi_\omega(a)\Omega_\omega\rangle \\
&=& \inf_{a \in A_0} \omega_n(|d + a|^2) \\
&=& \inf_{a \in A_0} \tau(h|d + a|^2).
\end{eqnarray*}

On taking the infimum over all $d \in \mathcal{D}$ with $\Delta(d)
\geq 1$, the result follows from Corollary \ref{vnc2}.
\end{proof}

{\bf Remark.}   After seeing this last result, Xu was able to use our
Szeg\"o formula, and facts about singular states, to remove the
hypothesis that $\mathrm{dim}(\mathcal{D}) < \infty$, and
to replace the `$2$' by a general $p$.  See \cite{BX} for details.

\section{Inner-outer factorization and the characterization of
outers}

Throughout this section $A$ is a maximal subdiagonal algebra. We
recall that if $h \in H^1$ then  $h$ is {\em outer} if $[h A]_1 =
H^1$.
An {\em inner} element is a unitary which happens to be in $A$.

\begin{lemma} \label{L3} Let $1 \leq p \leq \infty$.
Then $h \in L^p(M)$ and
 $h$ is outer in $H^1$, iff
$[h A]_p = H^p$.  (Note that $[\cdot]_\infty$ is the weak* closure.)

If these hold, then $h \notin [hA_0]_p$.  \end{lemma}

\begin{proof} It is
obvious that if $[h A]_p = H^p$ then $[h A]_1 = H^1$. Conversely, if
$[h A]_1 = H^1$ and $h \in L^p(M)$, then
the first part of the proof of \cite[Lemma 4.2]{BL3} applied to
$[hA]_p$ actually shows that $[h A]_p = [h A]_1 \cap L^p(M)$ for all
$1 \leq p \leq \infty$. Hence by \cite[Proposition 2]{Sai}, we have
$$[h A]_p = [h A]_1 \cap L^p(M) = H^1 \cap L^p(M) = H^p.$$

If $h \in [hA_0]_p$ then $1 \in [h A]_p \subset [[h A_0]_p A]_p
\subset [hA_0]_p$. Now $\Phi$ continuously extends to a map which
contractively maps $L^p(M)$ onto $L^p({\mathcal D})$ (see e.g.\
Proposition 3.9 of \cite{MW}).  If $h a_n \to 1$ in $L^p$, with $a_n
\in A_0$, then
$$\Phi(h a_n) = \Phi(h) \Phi(a_n) = 0
\to \Phi(1) = 1 ,$$
This forces $\Phi(1) = 0$, a contradiction.
\end{proof}

\begin{lemma} \label{L1}  If $h \in H^1$ is outer then
as an unbounded operator $h$ has dense range and trivial kernel.
Thus  $h = u |h|$ for a unitary $u \in M$.   Also, $\Phi(h)$ has
dense range and trivial kernel.
\end{lemma}

\begin{proof}  If $h$ is considered as an unbounded operator,
and if $p$ is the range projection of $h$, then since there exists a
sequence $(a_n)$ in $A$ with $h a_n \to 1$ in $L^1$-norm, we have
that $p^\perp = 0$. Thus the partial isometry $u$ in the polar
decomposition of $h$ is an isometry, and hence is a unitary, in $M$.
It follows that $|h|$ has dense range, and hence it, and $h$ also,
have trivial kernel.

For the last part note that $$L^1({\mathcal D}) = \Phi(H^1) =
\Phi([h A]_1) = [\Phi(h) {\mathcal D}]_1.$$ Thus we can apply the
above arguments to $\Phi(h)$ too.
 \end{proof}

There is a natural equivalence relation on outers:

\begin{proposition} \label{L2}  If $h \in H^p$ is outer,
and if $u$ is a unitary in ${\mathcal D}$, then $h' = u h$ is outer
in $H^p$ too.  If $h, k \in H^p$ are outer,  then $|h| = |k|$ iff
there is unitary $u \in {\mathcal D}$ with $h = u k$.  Such $u$ is
unique.
\end{proposition}

\begin{proof}  The first part is just as in the classical case.
If $h, k \in H^1$ and $|h| = |k|$, then it follows, as in \cite[p.\
228]{SW}, that $h = u k$ for a unitary $u \in M$ with $u, u^* \in
H^1$. Thus $u \in H^1 \cap M = A$ (see \cite{Sai}),
and similarly $u^* \in A$, and so $u \in {\mathcal D}$.
The uniqueness of $u$ follows since the left support projection of
an outer is $1$ (see proof of Lemma \ref{L1}).
\end{proof}

{\bf Remark.} \ It follows that if $f = uh$ is a
`Beurling-Nevanlinna factorization' of an $f \in L^p(M)$, for a
unitary $u \in M$ and $h$ outer in $H^1$, then this factorization is
unique up to a unitary in ${\mathcal D}$. For if $u_1 h_1 = u_2 h_2$
were two such factorizations, then $|h_1| = |u_1 h_1| = |u_2 h_2| =
|h_2|$. By Proposition \ref{L2} we have $h_1 = u h_2$ where $u \in
{\mathcal D}$ is unitary.
 Since $u_1 h_1 = u_1 u h_2  = u_2 h_2$, and
 since $h_2$ has dense range (see Lemma
\ref{L1}), we conclude that $u_2 = u_1 u$ and $h_2 = u^* h_1$.

\medskip

As in the classical case, if $h \in H^2$ is outer, then $h^2$ is
outer in $H^1$.  Indeed one may follow the proof on \cite[p.\
229]{SW}, and the same proof shows that a product of any two outers
 is outer (see also the last lines of the proof
of Theorem \ref{chr0} below).
We do not know whether every outer in $H^1$ is the square of an
outer in $H^2$.

The first main theorem of the present section is a generalization of
the beautiful classical characterization of outers in $H^p$:

\begin{theorem}  \label{chr0}  Let $A$ be a subdiagonal
algebra, let $1 \leq p \leq \infty$ and $h \in H^p$. If $h$ is outer then
$\Delta(h) = \Delta(\Phi(h))$. If $\Delta(h) > 0$, this condition is also
sufficient for $h$ to be outer.
\end{theorem}

Note that if $\mathrm{dim}({\mathcal D}) < \infty$, then  $\Phi(h)$
will be invertible for any outer $h$ by Lemma \ref{L1}. Thus in this
case it is automatic that $\Delta(\Phi(h)) > 0$.

\begin{proof} The case for general $p$ follows from the $p = 1$ case
and Lemma \ref{L3}. Hence we may suppose that $p = 1$.

First suppose that $h$ is outer. Given any $d \in L^1(\mathcal{D})$
and any $a_0 \in [A_0]_1$, we clearly have that $\tau(|d - a_0|)
\geq \tau(|d| - u^*a_0) = \tau(|d|)$, where $u$ is the partial
isometry in the polar decomposition of $d$. In other words, for any
$a \in [A]_1$ we have
$$\|\Phi(a)\|_1 = \inf_{a_0 \in A_0} \|a - a_0\|_1 =
\inf_{a_0 \in [A_0]_1} \|a - a_0\|_1 .$$ Therefore
$$\tau(|\Phi(h)\widetilde{d}|) = \inf_{a_0 \in A_0}
\tau(|h\widetilde{d} - a_0 |) , \qquad \widetilde{d} \in
\mathcal{D}.$$ Notice that $[hA_0]_1 = [[hA]_1A_0]_1 = [[A]_1A_0]_1
= [A_0]_1$. Thus the above equality may alternatively be written as
$$\tau(|\Phi(h)\widetilde{d}|) = \inf_{a_0 \in A_0}
\tau(|h\widetilde{d} - ha_0|) , \qquad \widetilde{d} \in
\mathcal{D}.$$Finally notice that $|\Phi(h)\widetilde{d}| =
||\Phi(h)|\widetilde{d}|$ and  $|h(\widetilde{d} - a_0)| =
||h|(\widetilde{d} - a_0)|$. Therefore if now we take the infimum
over all $\widetilde{d} \in \mathcal{D}$ with $\Delta(\widetilde{d})
\geq 1$, Szeg\"o's theorem will force $$\Delta(\Phi(h)) =
\Delta(|\Phi(h)|) = \Delta(|h|) = \Delta(h) .$$

Next suppose that $\Delta(h) = \Delta(\Phi(h)) > 0$. We will first
consider the case $h \in [A]_2$. By Lemma \ref{L3} we then need only
show that $h$ is outer with respect to $[A]_2$.  Replace $h$ by
$\widetilde{h} = u^*h$ where $u$ is the partial isometry in the
polar decomposition of $\Phi(h)$. If we can show that
$\widetilde{h}$ is outer, it (and hence also $u^*$) will have dense
range, which would force $u^*$ to be a unitary. Thus $h =
u\widetilde{h}$ would then also be outer. Now notice that by
construction $|h| \geq |\widetilde{h}|$ and $\Phi(\widetilde{h}) =
|\Phi(h)|$. From this and the generalized Jensen inequality we have
$$\Delta(h) = \Delta(|h|) \geq \Delta(|\widetilde{h}|) =
\Delta(\widetilde{h}) \geq \Delta(\Phi(\widetilde{h})) =
\Delta(\Phi(h)) = \Delta(h).$$Thus $\Delta(\widetilde{h}) =
\Delta(\Phi(\widetilde{h})) > 0$. We may therefore safely pass to
the case where $\Phi(h) \geq 0$. By multiplying with a scaling
constant, we may also clearly assume that $\Delta(h) = 1$.

For any $d \in \mathcal{D}$ and $a \in A_0$ we have
\begin{equation}\label{szout}
\tau(|1 - h(d+a)|^2) = \tau(1 - \Phi(h)d - d^*\Phi(h)) + \tau(|h(d+a)|^2).
\end{equation}
To see this, simply combine the fact that $\tau \circ \Phi = \tau$ with the
observation that $\Phi(h(d+a)) = \Phi(h)\Phi(d+a) = \Phi(h)d$. With $d, a$ as
above, notice that $\tau(|h(d+a)|^2) = \tau(||h|(d+a)|^2)$.  By Szeg\"o's
theorem in the form of Corollary \ref{vnc2}, we may select sequences $(d_n)
\subset \{ d \in \mathcal{D}^{-1} : \Delta(d) \geq 1\}$, $(a_n) \subset A_0$,
 such that
$$\lim_{n \to \infty}\tau(|h(d_n+a_n)|^2) = \Delta(|h|^2) =
\Delta(h)^2 = 1 .$$  Claim: we may assume the $d_n$'s to be positive. To see
this, notice that the invertibility of the $d_n$'s means that for each $n$
we can find a unitary $u_n \in \mathcal{D}$ so that $d_nu_n = |d_n^*|$.
Since for each $n$ we have
$$\tau(|h(d_n+a_n)|^2) = \tau(|h(d_n+a_n)u_n|^2) =
\tau(|h(|d_n^*|+a_nu_n)|^2) ,$$
 the claim follows. Notice that then $\tau(\Phi(h)d_n) =
\tau(d_n^{1/2}\Phi(h)d_n^{1/2}) \geq 0$. Using in turn the
$L^2$-contractivity of $\Phi$, the fact that $\Phi(h(d_n+a_n)) =
\Phi(h)d_n$, and H\"older's inequality, we conclude that
$$\tau(|h(d_n+a_n)|^2) \geq \tau(|\Phi(h)d_n|^2) \geq \tau(|\Phi(h)d_n|)^2 \geq
\tau(\Phi(h)d_n)^2 \geq \Delta(\Phi(h))^2 = 1.$$Since $\lim_{n \to
\infty}\tau(|h(d_n+a_n)|^2) = 1$, we must therefore also have
$\lim_{n \to \infty}\tau(\Phi(h)d_n) = 1$. But if this is the case
then equation (\ref{szout}) assures us that $h(d_n+a_n) \to 1$ in
$L^2$-norm as $n \to \infty$. That is, $1 \in [hA]_2$. Clearly $h$
must then be outer.

Now let $h \in [A]_1$. By noncommutative Riesz factorization (see \cite{MW}) we may select $h_1, h_2 \in [A]_2$ so that $h = h_1h_2$. Since $\Delta(h_1)\Delta(h_2) =
\Delta(h) = \Delta(\Phi(h)) = \Delta(\Phi(h_1))\Delta(\Phi(h_2)) > 0$ and
$\Delta(h_i) \geq \Delta(\Phi(h_i))$ for each $i = 1, 2$ (by the generalized
Jensen inequality), we must have $\Delta(h_i) = \Delta(\Phi(h_i))$ for each
$i = 1, 2$. Thus both $h_1$ and $h_2$ must be outer elements of $[A]_2$.
Consequently
$$[hA]_1 = [h_1h_2A]_1 = [h_1[h_2A]_2]_1 = [h_1[A]_2]_1 =
[[h_1[A]_2]_2]_1 = [[A]_2]_1 = [A]_1, $$
so that $h$ is outer as required.
\end{proof}

{\bf Remark.} In the general non-antisymmetric case, one can have
outers with $\Delta(h) = 0$. Indeed in the case that $A = M =
L^\infty[0,1]$, then outer functions in $L^2$ are exactly the ones
which are a.e.\ nonzero. One can easily find an increasing function
$h : [0,1] \to (0,1]$ satisfying $\Delta(h) = 0$, or equivalently
$\int_0^1 \, \log h = - \infty$.  See also \cite{Nak}.  This prompts
the following:

\medskip

{\bf Definition.}  We say that $h$ is strongly outer if it is outer
and $\Delta(h) > 0$.

\smallskip

Note that if $\mathrm{dim}({\mathcal D}) < \infty$, then every outer
$h$ is strongly outer.

\begin{corollary} Let $1 \leq p, q, r \leq \infty$ with $\frac{1}{p} =
\frac{1}{q} + \frac{1}{r}$ and let $h = h_1h_2$ with $h_1 \in H^q$ and
$h_2 \in H^r$. If $\Delta(h) > 0$ then $h$ is outer in $H^r$ iff both
 $h_1$ and $h_2$ are outer.
\end{corollary}

\begin{proof}    Clear from the last theorem, results in Section 2, and
the generalized Jensen inequality (as in the last paragraph of the
proof of the last theorem).
\end{proof}

By the Riesz factorization mentioned in the introduction, any
$h \in H^r$ is a product of the form in the last result.

\begin{corollary}  \label{Dout}
 If $f \in L^p({\mathcal D})$ with $\Delta(f) > 0$
then $f$ is outer.
Indeed there exist $d_n \in {\mathcal D}$ with
$\Delta(f d_n) \geq 1$, and  $f d_n \to 1$ in $2$-norm.
Also, any $f \in L^p(M)$ with $\Delta(f) > 0$
has left and right support projections equal to $1$.  That is, as an
unbounded operator it is one-to-one and has dense range.
\end{corollary}

\begin{proof}
For the first assertion note that $\Phi(f) = f$
and so $\Delta(f) = \Delta(\Phi(f)) > 0$.  An inspection of the
proof of the theorem gives the $d_n$ with the asserted properties.
Thus $f$ clearly has left support projection $1$, and by symmetry
the right projection is $1$ too.
Finally note that for the last assertion we may assume that $M = {\mathcal D}$.     \end{proof}

\begin{corollary} \label{opp}
If $1 \leq p \leq \infty$ and $\Delta(h) > 0$
then $h$ is outer in $H^p$ iff $[Ah]_p = H^p$.
 \end{corollary}

\begin{proof}
Replacing
$A$ by $A^*$, it is trivial to see that
$\Delta(h) = \Delta(\Phi(h)) > 0$, is equivalent to
$\Delta(h^*) = \Delta(\Phi(h^*)) > 0$.  The latter is equivalent to
$h^*$ being outer in $H^2(A^*) = (H^2)^*$; or equivalently,
 to
$(H^2)^* = [h^* A^*]_2$.  Taking adjoints again gives the result.
\end{proof}

{\bf Remark.} \ The last result has the consequence that the theory
has a left-right symmetry; for example our inner-outer
factorizations $f = uh$ below may instead be done with $f = h u$  (a
different $u, h$ of course).

\medskip

The following is a variant of Theorem \ref{chr0}:

\begin{proposition} \label{tona}   If $h \in H^2$, then $h$ is outer iff
 the wandering subspace of $[hA]_2$ (see {\rm \cite{N,BL5}}) has a
 separating cyclic vector for the ${\mathcal D}$ action,
and $$\Vert \Phi(h) \Vert^2 = \inf \{\tau(|h (1 - a_0)|^2) : a_0 \in
A_0 \} .$$
\end{proposition}

\begin{proof}   (Following \cite{Nak}.) \
For $x \in L^1(M)$ set $\delta(x) = \inf \{\tau(||x|^{\frac{1}{2}}
(1 - a_0)|^2) : a_0 \in A_0 \} .$

First suppose that $h \in H^2$ is outer. Then $[hA]_2 \ominus
[hA_0]_2 = H^2 \ominus [A_0]_2 = L^2({\mathcal D})$, which has a
separating cyclic vector.
 We next prove that if $h \in H^2$ is outer, then
$\Vert \Phi(h) \Vert^2 = \delta(|h|^2)$.   To do this
we view $\Phi$ as the orthogonal projection from $L^2(M)$ onto
$L^2({\mathcal D})$, which restricts to
 the orthogonal projection $P$ from $[A]_2$ onto
$L^2({\mathcal D})$.   For any orthogonal projection $P$ from a Hilbert space
onto a subspace $K$, we have $\Vert P(\zeta) \Vert =
\inf \{ \Vert \zeta - \eta \Vert :
\eta \in K^\perp \}$.  Thus $\Vert \Phi(h) \Vert^2 =
\inf \{\tau(|h - a_0|^2) : a_0 \in [A_0]_2 \}$.
Since $h$ is outer, we have
 $[[h A]_2 A_0]_2 = [H^2 A_0]_2$, or $[h A_0]_2 = [A_0]_2$.
 Thus $\Vert \Phi(h) \Vert^2 =
\inf \{\tau(|h - h a_0|^2 : a_0 \in A_0 \} = \delta(|h|^2)$.

Conversely, suppose that  the wandering subspace of
 $[hA]_2$ has a separating cyclic vector.
By  \cite[Corollary 1.2]{BL3}, we have $[hA]_2 = u H^2$ for a
unitary $u \in [hA]_2 \subset H^2$.  We have $h = u k$,
with $k \in H^2$, and $[k A]_2 = u^* [h A]_2 = H^2$.   So $k$ is
outer.
If $\Vert \Phi(h) \Vert^2 = \delta(h)$, then using the notation in the last
paragraph,
$$\Vert \Phi(u) \Phi(k) \Vert^2 = \delta(|u k|^2) = \delta(|k|^2)
= \Vert \Phi(k) \Vert^2.$$
That is, $\tau(\Phi(k)^* (1 - \Phi(u)^* \Phi(u) ) \Phi(k)) = 0$.
Since by Lemma \ref{L1} the left support projection
of $\Phi(k)$ is $1$, the
functional $a \to \tau(\Phi(k)^*a\Phi(k))$ is faithful
on $M$ (indeed,  $\tau(\Phi(k)^* a \Phi(k))
\neq 0$ for any non-zero $a \in M_+$), which forces $\Phi(u)^* \Phi(u) = 1$.
A simple computation shows that
$\Phi(|u - \Phi(u)|^2)  = 1 - \Phi(u)^* \Phi(u) = 0$, so that
$u = \Phi(u) \in {\mathcal D}$.   Thus $h = uk$ is outer.
\end{proof}

A classical theorem of Riesz-Szeg\"{o} states that if $f \in L^1$
with $f \geq 0$, then $\int \, \log f \, > - \infty$ iff $f = |h|$
for an outer $h \in H^1$ iff $f = |h|^2$ for an outer $h \in H^2$.
We now turn to this issue in the noncommutative case.

In what follows we are adapting ideas of Helson-Lowdenslager and
Hoffman:

\begin{lemma}  \label{Fact1}  Suppose that $A$ is a maximal
 subdiagonal algebra, and that $k \in L^2(M)$ with
$k \notin [kA_0]_2$.  Let $v$ be the orthogonal projection of $k$
onto $[kA_0]_2$. Then $|k-v|^2 = \Phi(|k-v|^2) \in L^1({\mathcal
D})$.  Also,  $\Delta(|k-v|) \geq \Delta(k)$.
 \end{lemma}

\begin{proof}  Suppose that $k a_n \to v$, with $a_n \in A_0$.
Clearly $k - v \perp k (1-a_n) a_0 \in k A_0$ for all $a_0 \in A_0$.
In the limit, $k - v \perp (k - v) a_0$.  That is, $\tau(|k-v|^2
a_0) = 0$, which by \cite[1.1 \& 4.1]{BL2} implies that $|k-v|^2 =
\Phi(|k-v|^2) \in L^1({\mathcal D})$.

For the last assertion, note that by Lemma \ref{vnc}
we have $\Delta(|k-v|^2) = \inf \{ \tau(|(k-v)d|^2) :  d \in
\mathcal{D}_+ \; \textrm{with} \; \Delta(d) \geq 1  \}$.
But since $vd \in [kA_0]_2$ for every $d \in \mathcal{D}$, we may apply
Szeg\"o's theorem to conclude that this infimum majorises
$$\inf \{ \tau(|kd -k a_0|^2) : d \in \mathcal{D}_+ \;
\textrm{with} \; \Delta(d) \geq 1 , a_0 \in A_0 \} = \Delta(|k|^2) =
\Delta(k)^2 ,$$ using the fact that $|kd -k a_0| = ||k|(d - a_0)|$.
  \end{proof}

\begin{theorem}  \label{Fact2}  Suppose that $A$ is a
maximal subdiagonal algebra, and that $k \in L^2(M)$.
Let $v$ be the orthogonal projection of $k$ onto $[kA_0]_2$.   If
$\Delta(k) > 0$, then $k$ has an (essentially unique)
Beurling-Nevanlinna factorization $k = u h$, where $u$ is a unitary
in $M$, and equals the partial isometry in the polar decomposition
of $k-v$, and $h$ is strongly outer and equals $u^* k$. We also have
$\Delta(k) = \Delta(k-v)$.  If
 $|k-v|$ is bounded below
then $(k - v) d = u$ for some $d \in {\mathcal D}$.
\end{theorem}

\begin{proof}   By  Corollary \ref{vnc3}, $k \notin [kA_0]_2$.
 By the Lemma,
$|k-v|^2  \in L^1({\mathcal D})$.   Let $u$ be the partial isometry
in the polar decomposition of $k-v$.  Since $\Delta(k-v) \geq
\Delta(k) > 0$ by the Lemma, we deduce from Corollary \ref{Dout}
that $u$ is surjective, and hence is a unitary.  In the case that
$|k-v|$ is bounded below let $d = |k-v|^{-1} \in {\mathcal D}_+$,
and then $u = (k-v) d$. Let $h = u^* k \in L^2(M)$.  We claim that
$\tau(u^* k a_0) = 0$ for all $a_0 \in A_0$, so that $h  = u^* k \in
L^2(M) \ominus [A_0^*]_2 = H^2$.   To see this, let $e_n$ be the
spectral projection of $|k-v|$ corresponding to the interval
$[0,1/n]$.  Then by elementary spectral theory, and since $k-v = u
|k-v|$, we have $1-e_n = |k-v| r$ for some $r \in {\mathcal D}$.
(Take $r = g(|k-v|)$ where $g$ is $\frac{1}{t}
\chi_{(\frac{1}{n},\infty)}$.)   Thus
$$\tau(a_0^* k^* u (1-e_n)) =  \tau(a_0^* k^* (k-v) r) = 0,$$
 since $k a_0 r^* \in [kA_0]_2$ and $k - v
\perp [kA_0]_2$.   On the other hand, by the Borel functional
calculus it is clear that $e_n \to e$ strongly, where $e$ is the
spectral projection of $|k-v|$ corresponding to $\{ 0 \}$. Since
$\Delta(|k-v|)  \geq \Delta(k) > 0$ by the Lemma, it is easy to see
by spectral theory
that $e = 0$ (this is essentially corresponds to
the fact that a positive function $f$ which is $0$ on a nonnull set
has $\int \log f = - \infty$).   We conclude that $\tau(a_0^* k^* u
e_n) \to 0$, and it follows that $$\tau(a_0^* k^* u) = \tau(a_0^*
k^* u e_n) + \tau(a_0^* k^* u (1-e_n)) = 0.$$

   To see that $u^* k$ is outer, we will use the
criterion in Theorem \ref{chr0}. We claim that $\Phi(u^* k) =
|k-v|$. To see this, note that by the last paragraph we have
$\tau(u^* x) = 0$ for any $x \in [kA_0]_2$, and in particular for $x
= v c$ for any $c \in {\mathcal D}$.   We have
$$\tau(\Phi(u^* k)c) = \tau(u^* k c) = \tau(u^* (k-v) c)
= \tau(|k-v| c) .$$ Since this holds for any $c \in {\mathcal D}$ we
have $\Phi(u^* k) = |k-v|$.   Thus we have
by the generalized Jensen inequality \ref{vncc2} that
$$\Delta(k) =  \Delta(u^* k) \geq \Delta(\Phi(u^* k))
= \Delta(|k-v|) \geq \Delta(k).$$
Hence $h = u^* k$ is outer by
Theorem \ref{chr0}.

The uniqueness follows from the remark after Proposition \ref{L2}.
  \end{proof}

\begin{corollary} \label{fact2}  Suppose that  $A$ is a
maximal subdiagonal algebra with ${\mathcal D}$ finite dimensional,
and that  $k \in L^2(M)$ with $\Delta(k)
> 0$.  Let $v$ be the orthogonal projection of $k$ onto
$[kA_0]_2$.   Then  $|k-v|$ is invertible, and all the conclusions
of the previous theorem hold.
\end{corollary}

\begin{proof}
  By the above, $|k-v| \in
L^1({\mathcal D}) = {\mathcal D}$, and  $\Delta(|k-v|) \geq
\Delta(k) > 0$.  Thus $|k-v|$ is invertible since $\mathcal{D}$ is
finite dimensional.
 The rest follows from the previous theorem.
\end{proof}

We next give a refinement of the `Riesz factorization'
into a product of two $H^2$ functions:

\begin{corollary} \label{ref} If $A$ is a
maximal subdiagonal algebra with ${\mathcal D}$ finite dimensional,
and if $f \in L^1(M)$ with $\Delta(f)
> 0$, then there exists an outer
$h_2 \in H^2$, an invertible $d \in {\mathcal D}$ with $\Delta(d) =
\frac{1}{\sqrt{\Delta(f)}}$, and an $h_1 \in [fA_0]_1$ such that
$f - h_1 \in L^2(M)$, and $f = (f - h_1) d h_2$. If also $f \in
H^1$, then this can be arranged with $h_1 \in H^1$, $\Phi(h_1) = 0$,
and $f - h_1 \in H^2$.
\end{corollary}

\begin{proof} Let $k =
|f|^{\frac{1}{2}}$.  By Corollary \ref{vnc3} we have $k \notin
[kA_0]_2$. If $u, v$ are as in Theorem \ref{Fact2}, and if  $f = w
|f| = w k^2$ is the polar decomposition of $f$, then $$f = (w  k  u)
(u^* k) = (w k (k - v)) d h_2 =  (f - h_1) d h_2$$ where $h_2 = u^*
k$ and $h_1 = w k v$.

If $k a_n \to v$ in $L^2$ norm, with $a_n \in A_0$, then $f a_n = w
k^2 a_n \to w k v$ in $L^1$ norm. Thus $h_1 \in [f A_0]_1$. Also,
$f - h_1 = w k u d^{-1}\in L^2(M)$ (recall that since $\mathcal{D}$ is
finite dimensional, $d^{-1} = |k - v| \in \mathcal{D}$). If $f \in H^1$,
then $h_1 \in [f A_0]_1 \subset H^1$, and $\Phi(h_1) = 0$.  So
$f-h_1 \in H^1 \cap L^2(M) \subset L^2(M) \ominus [A^*]_2 =  H^2$.
\end{proof}

\begin{corollary} \label{ref4} If $A$ is a
maximal subdiagonal algebra,
and if $f \in L^1(M)$ with $\Delta(f)
> 0$, then there exists a strongly outer $h \in H^1$, and a unitary
$u \in M$ with $f = u h$.
\end{corollary}

\begin{proof}  By the proof of Corollary \ref{ref}, and in
that notation, we have
$f = wku h_2$ for an outer $h_2$.  Note that $w$ is
a unitary, since $f$ has dense range (Corollary \ref{Dout}).
Since $\Delta(wku)
= \Delta(k) > 0$ (by facts in Section 2),
we have by
the last theorem that $wku = U h_1$ for a
unitary $U$ and strongly outer $h_1 \in H^2$. Let $h = h_1 h_2$.
\end{proof}

\begin{corollary} \label{pth2}  If
$A$ is a maximal subdiagonal algebra, and
 $f \in L^p(M)$ then $\Delta(f) > 0$ iff $f = u h$
for a unitary $u$ and a strongly outer $h \in H^p$. Moreover, this
factorization is unique up to a unitary in ${\mathcal D}$.
\end{corollary}

\begin{proof}  ($\Rightarrow$) \ By Corollary  \ref{ref4} we
obtain the factorization with outer $h \in H^1$. Since $|f| = |h|$
we have $h \in L^p(M) \cap H^1 = H^p$ (using \cite[Proposition
2]{Sai}), and $\Delta(h) > 0$.

($\Leftarrow$) \ We have $\Delta(f) = \Delta(u) \Delta(h) > 0$.

The uniqueness of the factorization was discussed after Proposition
\ref{L2}.
\end{proof}

{\bf Remark.}   The $u$ in the last result is necessarily in
$[fA]_p$ (indeed if $h a_n \to 1$ with $a_n \in A$, then $f a_n = u
h a_n \to u$).

\begin{corollary} \label{pth3}  If
$A$ is a maximal subdiagonal algebra, then  $f \in H^p$ with
$\Delta(f) > 0$ iff $f = u h$ for an inner $u$ and a strongly outer
$h \in H^p$. Moreover, this factorization is unique up to a unitary
in ${\mathcal D}$.   \end{corollary}

\begin{proof}  Clearly $f$ is also in $H^1$. Then $u$
is necessarily in $[fA]_p \subset H^1$, by the last remark. So $u \in H^1
\cap M = A$ (see \cite{Sai}).
 Thus $u$ is `inner' (i.e.\ is a unitary in $H^\infty = A$).
\end{proof}



An obvious question is whether there are larger classes of
subalgebras of $M$ besides subdiagonal algebras
for which such classical factorization theorems hold.
The following shows that, with a qualification, the answer to
this is in the negative:

\begin{proposition} \label{gian}  Suppose that $A$ is a tracial
subalgebra of $M$ in the sense of our previous papers, such that
every $f \in L^2(M)$ with $\Delta(f) > 0$ is a product $f = u h$
for a unitary $u$ and an outer $h \in [A]_2$.   Then
$A$ is a finite maximal subdiagonal algebra.
\end{proposition}  \begin{proof}
Suppose that $A$ is a tracial
subalgebra of $M$ with this factorization property.
  We will show
that $A$ satisfies the `$L^2$-density' and the
`unique normal state extension' properties which together
were shown in \cite{BL2} to imply that $A$ is subdiagonal.
As in our previous papers, $A_\infty$ is the tracial algebra
$A_\infty = M \cap [A]_2$ extending $A$.
If $x \in M$ is
strictly positive, then $\Delta(x) > 0$ by e.g.\ Theorem \ref{Fkd}
(2).  So $x = u h$ for a unitary $u$ and $h \in H^2$.  Clearly
$h$ is bounded, so that $h \in A_\infty$, and $x = (x^* x)^{\frac{1}{2}} =
|h|$.
Also, $h^{-1} \in A_\infty$, since
if $h a_n \to 1$ then $a_n \to h^{-1}$.
  Thus $A_\infty$ has the `factorization' property and so is
maximal subdiagonal \cite{BL2}.  Hence $A_\infty + A_\infty^*$, and therefore
also $A + A^*$, is dense in $L^2(M)$.  Next, suppose that
$g \in L^1(M)_+$ satisfies $\tau(g A_0) = 0$.
We need to show that $g \in L^1({\mathcal D})_+$. Since $\tau((g + 1) A_0) =
0$, we can replace $g$ with $g + 1$ if necessary, to ensure that
$\Delta(g) > 0$.
Let $f = g^{\frac{1}{2}} \in L^2(M)$.  Then $\Delta(f) > 0$,
$f \perp [fA_0]_2$,
and by hypothesis $f = u h$ for an outer  $h \in [A]_2$ and
some unitary $u$ in $M$.   Since $h = u^* f \perp u^*[fA_0]_2 =
[h A_0]_2 = [A_0]_2$, and $h \in [A]_2$,
it follows that $h \in [{\mathcal D}]_2$.   Thus
$g \in [{\mathcal D}]_1 = L^1({\mathcal D})$.   This verifies
the `unique normal state extension' property of \cite{BL2}.
\end{proof}

The following generalizes \cite[Theorem 5.9]{Ho}:

\begin{corollary} \label{pth}  If
$f \in L^1(M)_+$, then the following are equivalent:
\begin{itemize}
\item [(i)]  $\Delta(f) > 0$,
\item [(ii)] $f = |h|^p$ for a strongly  outer $h \in H^p$,
\item [(iii)]   $f = |k|^p$
for $k \in H^p$ with $\Delta(\Phi(k)) > 0$.
\end{itemize}
\end{corollary}

\begin{proof}  (i) $\Rightarrow$ (ii)  \  By
a previous result,  $\Delta(f^{\frac{1}{p}}) > 0$,
 and
so by the last result we have $f^{\frac{1}{p}} = u h$, where $h$ is
outer in $H^p$, and $u$ is unitary.  Thus $f = (f^{\frac{1}{p}}
f^{\frac{1}{p}})^{\frac{p}{2}} =
 (h^* h)^{\frac{p}{2}}  = |h|^p$.

(ii) $\Rightarrow$ (iii)  \ This follows from
Theorem \ref{chr0}.

(iii)    $\Rightarrow$ (i)  \ If $f = |k|^p$ for $k \in H^p$ with $\Delta(\Phi(k)) > 0$,
then $\Delta(f) = \Delta(k)^p \geq \Delta(\Phi(k))^p > 0$
by Theorem \ref{Fkd} and the generalized Jensen inequality.
\end{proof}

Of course in the case that ${\mathcal D}$ is finite dimensional
one can drop the word `strongly' in the last several results.
 In particular, in the case that the algebra $A$ is antisymmetric, these
results and their proofs
are much simpler and are spelled out in our survey \cite{BL5}.

\medskip

{\em Question.}  \ Is there a characterization of outers in $H^1$ in
terms of
extremals, as in the deLeeuw-Rudin
theorem of e.g.\ \cite[p.\ 139--142]{Hobk}, or \cite[pp.\ 137-139]{Gam}?

\medskip

{\bf Acknowledgements.}  We thank Barry Simon for helpful
information on Szeg\"o's theorem and its generalizations,
and for pointing out Verblunsky's precedence to the result usually
attributed to  Kolmogorov and Krein.
 We thank W. B. Arveson for continual encouragement,
Finally, we are grateful to Q. Xu for suggesting several years ago
to look at the noncommutative variant of the Riesz-Szeg\"o theorem;
and also for many insightful and valuable comments on the first
version of our paper. He together with Bekjan have recently
continued the $H^p$ theory contained in the present paper, by
extending it to values $0 < p < 1.$ This, together with other very
interesting related results of theirs, is contained in \cite{BX}.

\end{document}